\documentclass[twoside,11pt]{article} 

\usepackage{amssymb}
\usepackage{algorithm}
\usepackage{algpseudocode}
\usepackage{graphicx}
\graphicspath{{figures/}}
\newcommand{\ind}{\mathbf{1}}
\newcommand{\argmax}{\mbox{arg}\hspace{-2mm} \displaystyle\max} 

\title{A new gene co-expression network analysis based on Core Structure Detection (CSD)}

\author{A-C. Brunet($\dagger,\sharp$), J-M. Azais$(\sharp$), J-M. Loubes($\sharp$), J. Amar($\star,\ddagger$) and R. Burcelin($\ddagger$) 
\bigskip
\vspace{0.4cm} \\
\small $\dagger$ Vaiomer, Labege, France\\
\small $\sharp$ Institut de Math\'ematiques de Toulouse - IMT, Paul Sabatier University - Toulouse III, France\\
\small $\star$ F\'ed\'eration de cardiologie, P\^ole cardiovasculaire et m\'etabolique, CHU Toulouse, France\\
\small $\ddagger$ Institut National de la Sant\' et de la Recherche M\'edicale (INSERM), Toulouse, France,\\ 
\small Paul Sabatier University - Toulouse III,  (UMR \small 1048), Institut des Maladies M\'etaboliques et Cardiovasculaires (I2MC), \\
\small Team 2 : 'Intestinal Risk Factors, Diabetes, Dyslipidemia' }

\begin{document}

\maketitle

\begin{abstract}
We propose a novel method to cluster gene networks. Based on a dissimilarity  built using correlation structures, we consider  networks that connect all the genes based on the strength of their dissimilarity. The large number of genes require the use of the threshold to find sparse structures in the graph. in this work, using the notion of graph coreness, we identify clusters of genes which are central in the network. Then we  estimate a network that has these genes as main hubs. We use this new representation to identify biologically meaningful clusters, and to highlight the importance of the nodes that compose the core structures based on biological interpretations.   
\end{abstract}


\section{Background}

The analysis of hidden interaction structures in biological systems is  a problem of major importance. 
High-throughput technologies, and microarray data in particular, has allowed to obtain large quantities of gene expressions whose structure 
and relationships may explain the biological properties at stake.  
Interactions between genes are often described as a gene network, modeled by an undirected weighted graph.  
Each node or vertex represents a gene  while an edge connects two nodes if the corresponding genes are co-expressed, 
which means that their expression patterns are similar across the sample in the sense that the transcript levels increase or decrease in a similar way. 
Weights on the edges characterize the strength of the connections between the nodes. 
Rather than defining a similarity measure between each pair of variables without taking into account the others variables,  
the  graph structure adds information about adjacency, connectivity, and  characterize similarities within the complex structure of the gene network. 
The analysis of gene co-expression network made it possible to identify novel disease-related molecular mechanism as in~\cite{barabasi2004network}, \cite{weirauch2011gene}.
 \\
Identifying the most influential genes or centrally positioned genes in co-expression networks is the next problem.  
Indeed, the scale-free topological organization in biological networks is considered common property \cite{wang2003complex}, 
and is characterized by the presence of few nodes or hubs that  play a key role in the interaction mechanisms.  
These hubs  may play crucial roles in biological systems and may serve as valuable biomarkers for diseases see for instance 
in~\cite{he2006hubs}, \cite{jeong2001}. The hubs make sense as soon as the graph presents a structure of communities around these hubs,  
hence as soon as the graph can be partitioned into sub structures, each one containing a main hub. So clustering the graph into communities  
helps to summarize information by reducing the dimensionality from thousands of genes to a small number of cluster,
and in this way is an essential component of analysis methods, aiming at highlight the co-expression structure. 
Finally, clusters or communities of genes  exhibit similar expression patterns often related to functional complex \cite{eisen1998cluster},
as well as highly connected genes (hubs) in the network, may have critical functional roles or associated with key disease-related pathways. \\
Two problems are thus to be tackled : first the construction of the graph of the genes interactions and second the clustering of the network. 
\vskip .1in 
Building a graph to model the co-expression structure is not a trivial task. The starting point is the definition of a co-expression measure or similarity measure between the genes. The most popular measure is the absolute value of the Pearson's correlation coefficient \cite{allocco2004}, but  other measures can be considered.  A simple way to construct a graph based on pairwise similarities is to connect  all pairs of nodes having non-zero similarity values  and to weight edges by the corresponding similarity values. This approach leads to construct a weighted almost fully connected graph with few  real zero values in the similarity matrix, which prevents the discovery of particular clusters of interaction. Most of graph clustering methods are usable and efficient if the graph is sufficiently sparse,  such that the organization in clusters is rather pronounced with many edges within each cluster   and relatively few between the clusters. Hence, methods that promote the sparsity of the network have received a growing attention over the last decade. One way to estimate a sparse graph is to connect two nodes if their similarity exceeds a certain threshold $\epsilon$, defining an $\epsilon$-neighorhood graph.  Another way is to connect each node with all its k-nearest neighbors,  a $k$-nearest neighbor graph. Both approaches are very sensitive to the choice of the parameter $\epsilon$ or k. The advantage of a $k$-nearest neighbor graph, compared to an $\epsilon$-neighborhood graph, is  that it can connect nodes using different local similarity thresholds for each nodes 
(edge connecting nodes in each neighborhood can have highly variable weights), rather than connecting nodes based on a unique threshold. Yet the $\epsilon$-neighborhood graph estimation fails to identify regions or communities of a graph with different levels of density. For a given stringent threshold, many nodes in the weakly dense regions of the graph may be disconnected, while for a low threshold, nodes in the dense regions of the graph may have too many edges, making clustering in these regions difficult. Moreover, the $k$-nearest neighbor graph also have  a tendency to include non significant edges due to the use of a fixed-size neighborhood.\\ Identifying clusters of similar nodes of graph is also a hard task. The mere concept of graph clustering hides very difficult issues. 
As a matter of fact, clustering implies splitting the data into groups that share a common behavior. 
Hence it requires being able to define a notion of similarity or equivalently distance between the different elements.
When dealing with non Euclidean structures such as graphs for instance, this notion of similarity is not natural 
and highly relies on the objective of the study.  Actually several authors have proposed different ways of achieving this goal. 
We refer to\cite{fortunato2010community}, \cite{newman2004finding,newman2004detecting,newman2004fast} or \cite{schaeffer2007graph} for instance and references therein. 
A very common method employed for clustering biological networks is the hierarchical clustering.
In this class of methods we can cite the WGCNA algorithm \cite{langfelder2008wgcna,zhang2005general} 
widely used for the analysis of co-expression networks \cite{dileo2011weighted,ivliev2010coexpression}.
The starting point of hierarchical graph clustering is the definition of a similarity measure between graph objects or nodes, 
e.g. the topological overlap measure \cite{yip2007gene}.
With the agglomerative (bottom-up) strategy, each node is considered as a separate cluster (singleton) 
in the starting point of the clustering procedure, 
and then clusters are iteratively merged if their similarity is sufficiently high. 
Since clusters are merged based on their mutual similarity, a measure that estimates how two clusters are closed must be defined,  for example by computing the average similarity between nodes in one group and nodes in the other (average linkage clustering).   The procedure ends up with a unique cluster that contains all the vertices of the graph.
The partition strategy is top-down as it follows the opposite direction, but has been rarely used in practice.
Hierarchical clustering has the advantage that it does not require any assumptions on the number and the size of the clusters.
However, the results of the method may strongly depend on the similarity measure adopted,  while the graph partition relies on the choice of the cut height in the dendrogram which may lead to very different clusters, as in~\cite{langfelder2008defining}.
Spectral clustering \cite{ng2002spectral} \cite{von2007tutorial} 
is another very popular graph clustering method used in a variety of applications  \cite{ruan2007efficient}.
This method consists in transforming the representation of the graph to make the cluster properties of the initial data sets much more evident.
In particular, by replacing the weighted adjacency matrix of the graph 
with a subset of eigenvectors of the Laplacian matrix, a different representation of the graph is induced and  
can be useful to extract information on the cluster structure. 
The problem of clustering nodes is then reduced to a classical clustering problem defined in an Euclidean space.
The set of points in the new space (or elements of eigenvectors) can be clustered via standard clustering algorithm, like the $k$-means clustering.  
Finally a variety of criteria and methods have been proposed to allow automatic selection of the number of clusters \cite{tibshirani2001estimating}, \cite{white2005spectral}. 
However, the more noisy or ambiguous the clusters are, the less effective are these methods. 
In many cases, determining the optimal value for the number of clusters is a nontrivial or even impossible problem.\\ \indent 
Hence graph clustering can be achieved in several ways having several advantages and drawbacks. 
In the analysis of networks, the issue is even more difficult since the induced graph is unknown and has to be estimated, resulting to different possible graph partitions. 
Since for real data set, there is no ideal underlying graph, the purpose of the clustering procedure and the notion of the estimation of the graph are deeply linked.
\vskip .1in
In this paper we propose a new method, called CSD (for Core Structure Detection), to simultaneously select clusters of co-expressed genes and provide a sparse network of genes. For this we will identify clusters of genes which are strongly connected together which  will be called cluster core structures. Then we will build the graph onto such central hubs.  Hence the CSD algorithm falls into two steps. The first step identifies sets of strongly connected nodes that make good candidates to define the cores of clusters. 
The second step aggregates  the remaining nodes out of the core structures to build the final clusters. The algorithm takes as input the similarity matrix between the variables (e.g. matrix of absolute correlation coefficients for co-expression network)  and a parameter that characterizes the minimum size of these core structures. More precisely, we consider all the possible graphs obtained by removing stepwise the edges with the lowest weight and  identifying the set of nodes which  are preferentially connected on this collection of graphs. They will be denoted as the core structures. Then the nodes are connected to these core structures to obtain the clusters of the graph, with each cluster containing a unique core structure and all the nodes that are more strongly connected with nodes in this core structure than with nodes in a different core structure. This process allows the detection of core structures of varying sizes and densities governed by an input parameter that controls the level of granularity of the result partition. Note that the clusters obtained with the CSD algorithm, for different values of the input parameter or the minimum size of the core structures, are nested according to a hierarchical organization. A study of the evolution of the partition, as a function of the input parameter, enables to select  the best compromise with respect to the number and the size of the obtained clusters.
\vskip .1in
In the Results and Discussion section, we introduce the CSD approach before analyzing the performances of the algorithm,
and compare its performances with the WGCNA and spectral clustering algorithms.
For this purpose, we first use simulated data for explicit estimations of the performances (a priori knowledge of the community structure), and second, we use well-studied transcriptomic data describing gene expression 
during the cell division cycle of the yeast Saccharomyces cerevisiae \cite{spellman1998comprehensive}
to test the ability of the method to identify biologically meaningful clusters, and to highlight the importance of the nodes that compose the core structures based on biological interpretations.   In the following Methods section, we provide more details on the mathematical formalism of the network clustering algorithms used in the analysis and on the criteria used to analyze and compare the performances of network clustering algorithms.

\section{Results and Discussion}

\subsection{The CSD algorithm}

\subsubsection{Graph definition}
Our starting point is a set of joint expression, or co-expression, of $p$ genes over $N$ individuals.
From this co-expression network, a weighted adjacency matrix $W=(w_{ij})_{1\leq i,j \leq p}$ is constructed, 
in most of the case using the absolute value of the Pearson correlation coefficient.  Hence a gene co-expression network can be modeled as a weighted undirected graph in which each node corresponds to a gene 
while the edges describe gene co-expression relationships. $W$ is thus the  weighted adjacency matrix of a graph $G=(V,E)$ with  $V$ a set of $p$ nodes and  $E$  bidirectional edges. We assume in addition that no edge connects a node to itself. This assumption models the fact that we do not want a gene to   have interactions with itself.\vskip .1in
The weighted adjacency matrix $W$  gives the connection strength between each pair of nodes. \\
This network may exhibit too many connexions to be used by biologists, hence  we need to estimate a sparse adjacency matrix $W$  to obtain a network with fewer links between the genes. A natural procedure is to threshold the elements in the similarity matrix according to a given threshold parameter $\lambda$, using a hard thresholding scheme.
The new adjacency matrix $W_\lambda$ is then defined by 
$$
w_{ij}^\lambda = w_{ij}\ind_{w_{ij}\geq\lambda} \ 1\leq i,j\leq p
$$
The graph associated to $W_0=(w_{ij}^0)$ is fully connected while as $\lambda$ tends to $1$, 
the graph $G_\lambda$ associated to $W_\lambda$ consists of singletons. This procedure amounts to select a chosen number of correlations and deleting the smallest ones.

\subsubsection{What is a core structure ?}

In previous works, the main issue consists in choosing a good threshold parameter to obtain a unique graph $G_\lambda$. Here,
we emphasize that that genes in different functional pathways may exhibit different levels of co-expression, 
involving the formation of subgraphs with different levels of density or similarities. Hence, rather than  constructing a single sparse graph to highlight the modular structure, we propose to study the evolution of the sparsity of the graphs defined by removing step by step edge with lowest weight,  which amounts to threshold the similarity matrix with a growing parameter.  Then we identify the groups of vertices that stay connected between themselves while they are disconnected in the graphs obtained for  larger  threshold values. \vskip .1in
We first introduce the notion of central group before defining the core structure notion.
We call ``n-central group'' a set containing at least n nodes, which are  strongly connected between themselves  than to other nodes out of the set, hence more similar with the nodes in their group. So there exists a graph $G_\lambda$ where this group of nodes forms a connected component  
(any two vertices are connected to each other by a path in the graph and there is no connection with the vertices out of the component). This notion is similar to the notion of coreness  introduced  in the context of social networks in~\cite{MR721295}.
\paragraph{Definition 1} Let $W=(w_{ij})_{1\leq i,j \leq p}$ be a weighted adjacency matrix such that $w_{ij}\in [0,1]$, $\forall i,j \in \{1,2,...,p\}$.
For any $n\in \mathbb{N}^*$, a n-central group $C$ is a subset of elements, 
$C\subset \{1,2,...,p\}$, satisfying the following conditions
\begin{enumerate}
\item[(i)] $|C|\geq n$ 
\item[(ii)] There exists $\lambda\in [0,1)$ such that $C$ is a connected component of the graph $G_\lambda$.
\end{enumerate}
Then we define a ``n-core structure'' as a maximal subset of elements 
that have the property to be a n-central group which can not be split into two or more disjoint n-central groups.
\paragraph{Definition 2}  Let $W=(w_{ij})_{1\leq i,j \leq p}$ be a weighted adjacency matrix such as $w_{ij}\in [0,1]$, $\forall i,j \in \{1,2,...,p\}$.
For any $n\in \mathbb{N}^*$, a n-core structure $Q$ is a maximal subset of elements, 
$Q\subset \{1,2,...,p\}$ satisfying the following conditions
\begin{enumerate}
\item[(i)] $Q$ is a n-central group  
\item[(ii)] For all $\lambda\in [0,1]$, the elements of $Q$ cannot be split into two or more connected components of $G_\lambda$ 
containing at least $n$ nodes (n-central group). 
\end{enumerate}

Intuitively, a core structure is a good candidate to be centrally positioned in a cluster of co-expressed genes.
It appears, indeed, that the core structures group the most significantly co-expressed genes at all levels of co-expression.
\\ \newline
The algorithm for detecting the core structures (presented in the Methods section) just perform hierarchical clustering 
by removing edges between pairs of vertices with low similarities and then
analyze the evolution of the structures get disconnected from each other. \\
\indent At the beginning  of the algorithm, we start with a graph $G_{\lambda_0}$ where $\lambda_0=0$ which is the fully connected graph ($W_{\lambda_0}$).  At each step $t$, a new graph $G_{\lambda_t}$ is defined by removing the edge with lowest weight on $G_{\lambda_{t-1}}$ 
i.e. $\lambda_t=\displaystyle \mbox{min}_{\{i,j | w_{ij}^{t-1}\neq 0 \}} w_{ij}^{t-1}$.
It is not necessary to study all possible graphs $G_\lambda$ obtained with 
$\lambda$ in the set of nonzero zero values of the similarity matrix S.
We will only consider paths of maximum capacity or paths maximizing the weight of the minimum edge in the path. So two nodes are in the same connected component of the graph $G_\lambda$ 
as long as the threshold parameter $\lambda$ is lower than the capacity of the maximum capacity path between these nodes.
The graph $G_t$ has one more connected component  than the graph $G_{t-1}$ if $\lambda_t$ is equal to the minimum capacity 
of all maximum capacity paths in $G_{t-1}$.
We will then define the maximum spanning tree (each path between two nodes is the maximum capacity path between these node) 
of the graph $G_{\lambda_0}=G_0$
and redefine the weighted adjacency matrix $W_{\lambda_0}$ such that 
$w_{ij}^{\lambda_0}=0$ if there is no edge between nodes i and j on the maximum spanning tree.
By this process we improve the running time of the algorithm for core structure detection from  $O(p^3)$ to $O(p^2)$.

\subsubsection{Clustering based on core structures}

Once the n-core structures are detected as previously described, we use it to initialize cluster centers and 
we propose a simple stepwise algorithm for reconnecting nodes around the core structures and then completing clusters. 
Suppose we have identify $K$ core structures $Q_1,.Q_2,..,Q_K$.
Each cluster is at first defined by a unique core structure,  $C_{1,0}=Q_1,C_{2,0}=Q_2,...,C_{K,0}=Q_K$.
At each iteration t, we denote by $U_t$ the set of unclustered element i.e. 
$i\in U_t\ \Rightarrow \ i\not\in C_{k,t-1},\forall k=1,2,...,K$, 
and we identify the unclustered element $i\in U_t$ and the clustered element $j\not\in U_t$ 
that have the maximum similarity value i.e. 
$s_{ij}\geq s_{i^\prime j^\prime},\ \forall\ i^\prime\in U_t,\ \forall\  j^\prime\not\in U_t$.
We then put the unclustered element $i$ in the cluster $C_{k,t-1}$ including the element $i$,
$C_{k,t}=C_{k,t-1}\cup \{i\}$ (and $C_{l,t}=C_{l,t-1},\ \forall l\neq k$).
The process is repeated until all elements are clustered and ensures that 
each node and its first nearest neighbor are necessarily placed in the same cluster.

\subsubsection{Parameter sensitivity analysis}

\begin{figure}[!h]
\centering
\begin{minipage}[b]{0.9\linewidth}
\centering
  \includegraphics[scale=0.6]{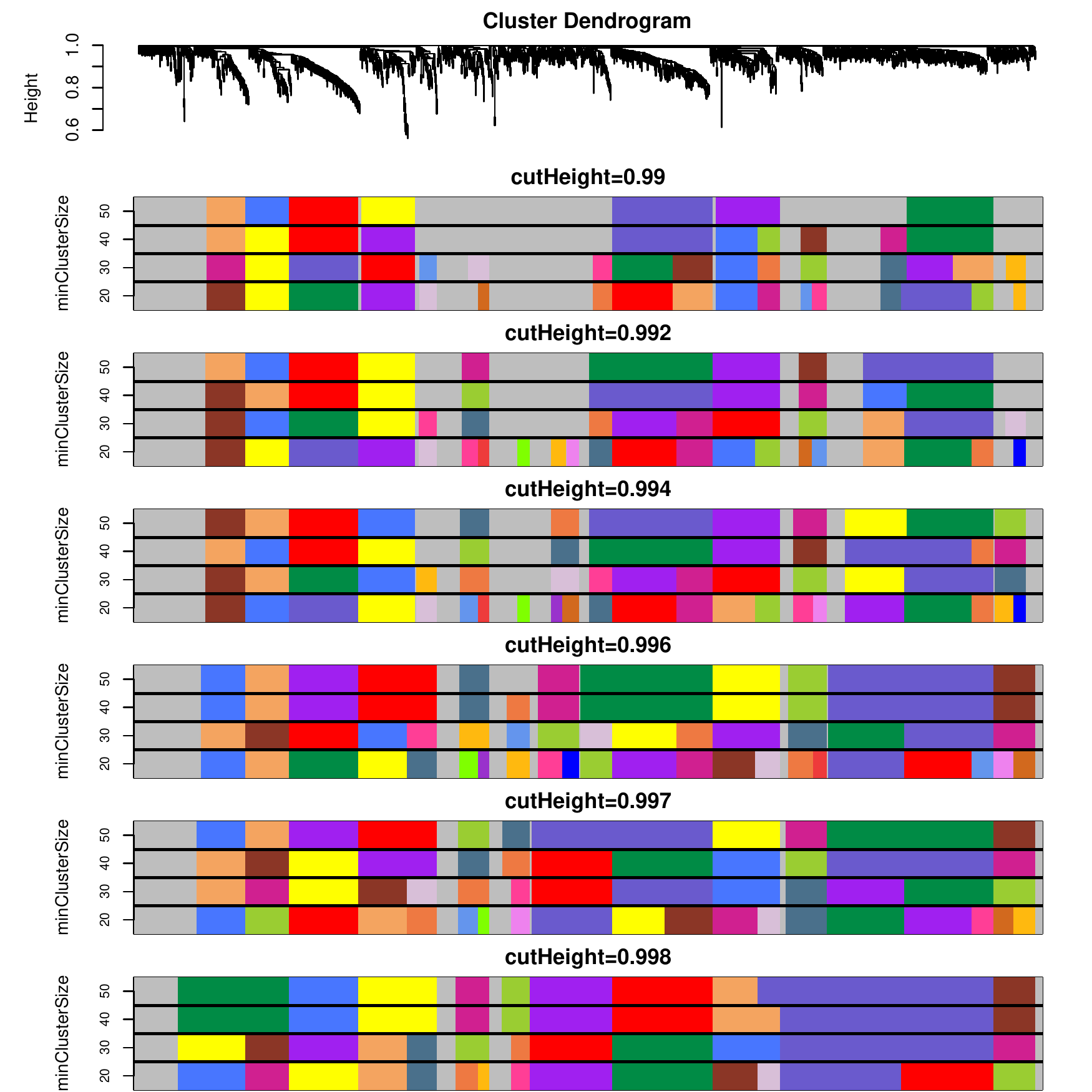}
  \caption{\label{wgcna_choose_para} \small WGCNA algorithm  applied to yeast cell-cycle gene expression data set.
  Representation of the hierarchical cluster tree and representation of the cluster results 
  for various dynamic tree cut input parameters (height cut and minimum size of clusters),  
  which are shown in the color bands.
  The gray color is reserved for non-clustered genes.}
  \end{minipage}
\end{figure}

\begin{figure}[!h]
\centering
\begin{minipage}[b]{0.8\linewidth}
\centering
  \includegraphics[scale=0.5]{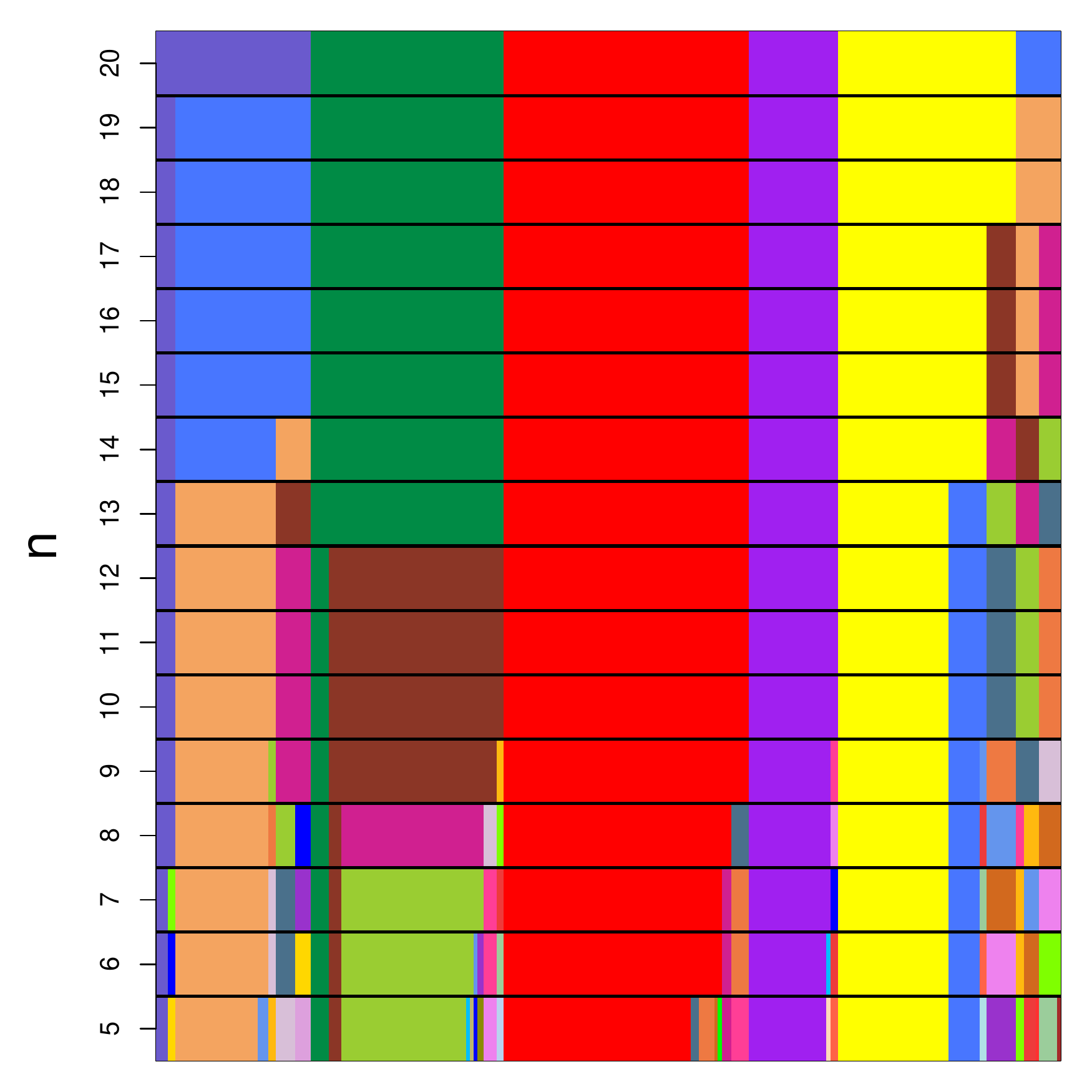}
  \caption{\label{core_choose_para} \small CSD algorithm applied to yeast cell-cycle gene expression data set 
  for various values of the minimum size of the core structures. 
  The color bands on the x-axis shows the cluster membership according to different parameter ($n$) settings on the y-axis.}
  \end{minipage}
\end{figure}

The greatest difficulty in the field of data clustering is the need for input parameters 
which can greatly influence the result.
In many applications the optimal values of these parameters are very difficult to determine. 
To illustrate the robustness of our CSD algorithm, we analysed the effect that changes in parameter settings have 
on clustering results compared to the input sensitivity of spectral and WGCNA algorithm.  
\\ \newline
The CSD and spectral algorithms need a unique tuning parameter while WGCNA algorithm needs more parameters.
Indeed, the choice of a threshold parameter is first required to define the adjacency matrix and then, 
 more than one other parameters are required to detect clusters in the dendrogram with a cut method
 less inflexible and more adapted than the use of a constant height cutoff value.
For both CSD and WGCNA algorithm, the set of clusters returned are hierarchically structured 
and the input parameters control the granularity of the clustering results. 
If the granularity is too fine, there will be too many small clusters with a high correlation among their elements,
and conversely, a too low granularity leads to product a few large clusters with a low correlation among their elements.
A well-suited parameter setting may be recognizable by visual inspection of the evolution of the level of granularity 
at different input parameters.  
The height and the minimum size of cluster parameters of the dynamic tree cut method (WGNCA) 
provide improved flexibility for tree cutting but it is not easy to find optimal cutting parameters (Figure~\ref{wgcna_choose_para}). 
The CSD algorithm leads to embedded clusters and the number of clusters is monotonically decreasing as the parameter $n$ increases, 
so clustering results with well-fitted levels of granularity can be quickly and easily identified (Figure~\ref{core_choose_para}).
The problem of the right choice of the number clusters for spectral clustering algorithm is crucial and still open issue.
Even if a variety of methods have been devised for this problem, 
most of this methods works well if the data contains very well pronounced clusters, 
and results are more ambiguous in cases of noisy or overlapping clusters. 
It should also be noted that the spectral clustering result using k-means algorithm depends largely on the initial set of centroids, 
so many iterations using random initial values are needed to have good performance (we used 500 random initial values).
\vspace{1cm}

The choice of the tuning parameter of the CSD algorithm is thus the one that promotes stability of the clusters as shown in Figure~\ref{core_choose_para}. This choice can easily be turned into an automatic one. Hence promoting coreness of cluster centers may be a good way to get stable structures for gene networks such as hubs and clusters. The following of the paper studies the biological consequences of such representation.


\subsection{Performance evaluation using simulated data}

To evaluate the ability of our algorithm to identify clusters or communities that are highly intra-connected 
in correlation graphs, we simulated correlated data such as the correlation matrix displaying community structure, 
with communities of varying sizes, densities, and level of noise. 
We also compared the accuracy and robustness of our algorithm with two very popular clustering algorithms, 
the spectral clustering algorithm and the WGCNA algorithm successfully applied in genes network clustering context. 

\subsubsection{A model to simulate co-expression matrix}
The goal of the analysis is to study the performances of our algorithm in many different situations, 
so we need to generate correlation graphs may contain communities or clusters of nodes with various size, 
intra and inter densities and levels of noise.  
A co-expression gene network may be modeled by a graph which is defined using microarray expression profiles, 
that are transformed to a co-expression matrix by calculating pairwise correlations. 
It is more convenient to simulate expression data instead of co-expression or correlation matrix directly, 
to create various source of randomness 
and obtain realistic co-expression matrices. 
We propose a model inspired by~\cite{langfelder2007eigengene} for simulating expression data sourced from different communities of genes, 
in such a way to define graphs resemble real gene co-expression network. 
Many complex networks and in particular, gene co-expression networks, are known to be scale free, meaning that the network has few highly 
connected nodes (hubs) and many poorly connected nodes. 
Simulation can be seen as a process through which we first generate expression profiles of 
leader genes or hubs, one hub for each community, 
and then, expression profiles of other genes are generated in each community on the basis of the corresponding hub, 
in order to have more or less high correlation with the hub 
(few high correlations and many moderate or small correlations in each community). 
Gene expression data is often assumed to be (log)normally distributed and 
the absolute value of the Pearson correlation coefficient is commonly used to measure the co-expression between genes. 

The aim is to generate gene expression profiles where each profile is a vector of size $N$.
Strategy for simulating gene profiles in a single community of size $n_C$ is as following

\begin{enumerate}
  \item Generate the expression profile $\mathbf{x}^{(1)}=(x_1^{(1)},x_2^{(1)},...,x_N^{(1)})^\prime$ 
 of one gene from a normal distribution, $x_i^{(1)}\sim \mathcal{N}(0,1),\forall i=1,2,...,N$. 
 This profile defines the leader or hub profile.
 \item Choose a minimum correlation $r_{min}$ and a maximum correlation $r_{max}$ between the leader gene and other genes in the community. 
 \item Generate $n_C-1$ expression profiles such that the correlation of the j-th profile $\mathbf{x}^{(j)}$ 
 with the  leader gene profile $\mathbf{x}^{(1)}$ is forall $j=2,3,...,n_C$ close to
 $$
 r_{j}=r_{min}+\Delta_r\left( 1-\frac{j}{n_C} \right),
 $$
 where $\Delta_r=r_{max}-r_{min}$.
 The $r_2$ correlation is close to $r_{max}$ and the $n_C$-th correlation is equal to $r_{min}+\Delta_r=r_{max}$.

 The j-th profile is generated by adding a Gaussian noise term to the leader profile,
 to allow a correlation with the leader profile close to the required correlation $r_j$
 $$
  x_i^{(j)}= x_i^{(1)}+\sqrt{\left( \frac{1}{r_j^2}-1 \right)} \epsilon_i^{(j)},\forall i=1,2,...,N
 $$
 where $ \epsilon_i^{(j)}\sim\mathcal{N}(0,1) $. 
\end{enumerate}

We created six particular scenarios for simulated co-expression data. 
For each scenario we fixed the size of profile $N=100$ and 
we simulated gene profiles within $K=5$ different communities of sizes $n_{C_1}$,$n_{C_2}$,...,$n_{C_K}$ 
randomly drawn from the set $\{50,100\}$. 
A simulated expression data set $X$ is then composed of $p=n_{C1}+n_{C_2}+...+n_{C_K}$ gene profiles
generated for the five communities with the process described above
$$X=\left(\mathbf{x}^{(1)}_{C_1},...,\mathbf{x}^{(n_{C_1})}_{C_1},\mathbf{x}^{(1)}_{C_2},...,\mathbf{x}^{(n_{C_2})}_{C_2},...,
\mathbf{x}^{(1)}_{C_K},...,\mathbf{x}^{(n_{C_K})}_{C_K}\right)\in\mathbb{R}^{N\times p}.$$
The coefficients $s_{ij}\in[0,1]$ of the $p\times p$ co-expression matrix $S$ are defined as the absolute values of 
Pearson correlation coefficients between gene profiles.

In a first scenario (S-1), we generated well-defined communities inside which correlations between gene profiles 
decrease linearly from $r_{max}=1$ to $r_{min}=0.5$. 
For the second scenario (S-2), we have allowed communities to have different levels of density
by randomly selecting, for each community, the parameter configuration  
$\left\{ r_{min}=0.5,\ r_{max} =1\right\}$ or the other configuration $\left\{ r_{min}=0.4,\ r_{max} =0.7\right\}$ (smaller density).
The third scenario (S-3) products two communities with high intercorrelation. 
We generated two leader profiles which have a correlation value close to $0.8$ (the three other profiles are independently drawn).
Genes inside communities are generated with the parameter configuration  $\left\{ r_{min}=0.5,\ r_{max} =1\right\}$.
The fourth scenario (S-4) consists in creating noisy data. We simulated expression data using  $\left\{ r_{min}=0.5,\ r_{max}=1\right\}$, 
then we have standardized the data and added to each gene profile a gaussian  noise term drawn from a normal distribution  $\mathcal{N}(0,1) $.
Another scenario (S-5) was designed to study the effects of the presence of irrelevant variables. 
A total of $p$ gene profiles are generated with parameters   $\left\{ r_{min}=0.5,\ r_{max}=1\right\}$ and 
$p$ irrelevant variables of size $N=100$ , drawn from a $\mathcal{N}(0,1)$, are added.   
The last scenario (S-6) is a mix of the five previous scenarios. 
We generated two leader profiles which have a correlation value close to r which is drawn from the set $\{0.2,0.4,0.6\}$
 and the three other profiles are independently drawn.
For each community, the parameter configuration is randomly selecting from the two configurations 
$\left\{ r_{min}=0.5,\ r_{max}=1\right\}$ and  $\left\{ r_{min}=0.4,\ r_{max}=0.7\right\}$.
Simulated data have been stantardized and then a gaussian noise were added to each  profile with a 
standard deviation drawn from the set $\{0.1,0.5,1\}$.
Additionnaly we added $p$ irrelevant profiles simulated from a $\mathcal{N}(0,1)$.

\begin{figure}[!h]
\centering
  \includegraphics[scale=0.58]{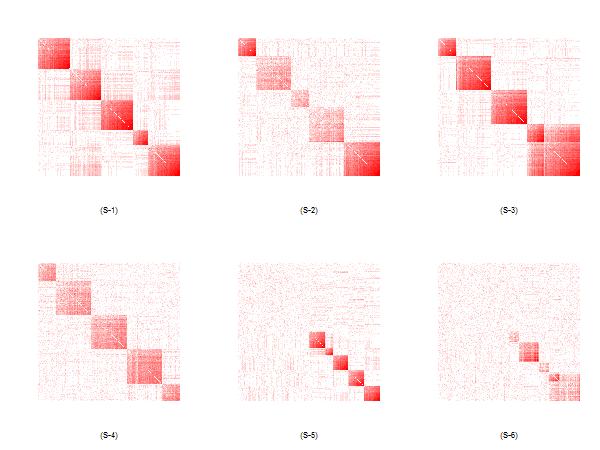}
  \caption{\label{simu_example} \small  Example of a co-expression matrix produced for the six scenarios.}
\end{figure}

We generated $100$ expression data sets for each scenario, and calculated corresponding co-expression matrices $S$ (absolute values of Pearson correlations).
Examples of co-expression matrices obtained for the six scenarios are shown in ~\ref{simu_example}.

\subsubsection{Results}

We studied the performance of the three clustering algorithms (CSD, spectral and WGCNA) on the 600 simulated co-expression matrices.
We needed to define a graph and choose appropriate parameter settings for each clustering algorithm. 
For the spectral clustering and the CSD algorithm we defined a fully connected graph by taking the weighted adjacency matrix of the graph 
equal to the co-expression matrix, W=S. 
For the use of  WGCNA algorithm we defined the weighted network adjacency matrix by 
raising the absolute value of the correlation matrix to a power, $w_{ij}=s_{ij}^\beta$,
 and chose the threshold parameter $\beta$ by applying the approximate scale-free topology criterion \cite{langfelder2008wgcna}.
Spectral clustering was run with $K=5$ (known parameter), and was repeated 50 times for each data set to avoid local minima.
More realistically, the number of classes is unknown, so we tried to find it 
by running the algorithm for a wide range of number of classes (from 2 to 10) 
and chose clustering result with the highest value of the Dunn cluster validity index \cite{dunn1973fuzzy}.
The WGCNA algorithm identifies classes using hierarchical clustering 
and the Dynamic cut algorithm \cite{langfelder2008defining}, 
is proposed for automatically detecting clusters in a dendrogram based on shape of its branches. 
We chose the minimum size of a community equal to $40$ and the default value  (0.99)  for the cut height
excepted for the  scenario S-4 for which it is better to choose a more high cut height value (0.999). 
Clustering with the CSD algorithm require to select a unique tuning parameter, the minimum core size, 
and we chose it equal to 10.\\

\begin{figure}[!h]
\centering
\begin{minipage}[b]{1\linewidth}
\centering
  \includegraphics[scale=0.5]{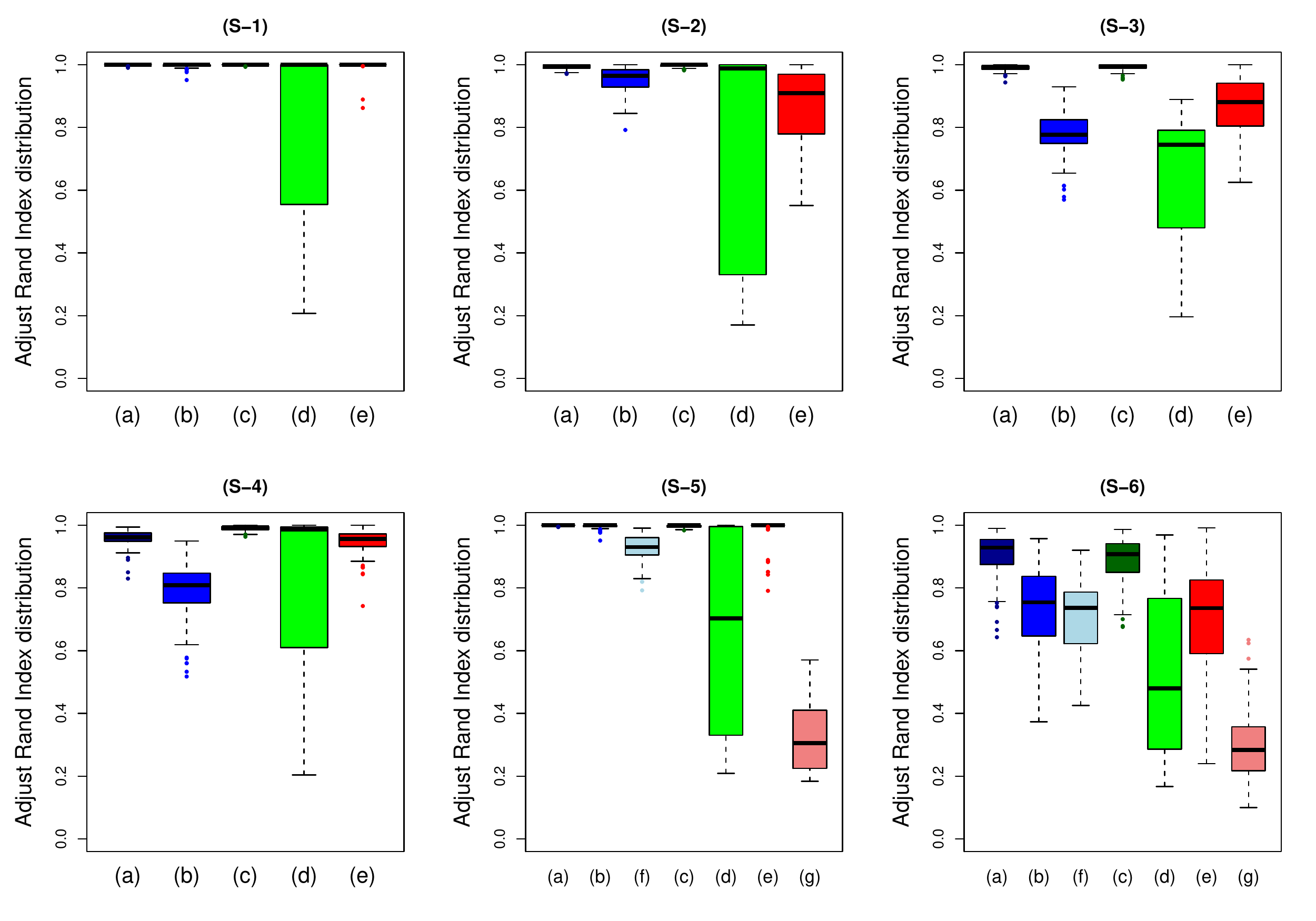}
   \caption{\label{simu_RI} \small Results on simulated data. 
  Evaluation of the clustering methods (CSD, spectral clustering and WGCNA) on the simulated datasets obtained with the six scenarios.
  The (a) boxes report adjust rand index for the CSD clustering results. 
 The (b) boxes report adjust rand index for the core structure labels identify with the CSD algorithm 
(one cluster is composed of the outer genes of the core structures).
  Evaluations of spectral clustering results correspond to the (c) boxes when the number of clusters is known,
and to the (d) boxes when the number of clusters is chosen based on the Dunn's index. 
The (e) boxes are reserved for WGNCA result evaluations. 
For the S-5 and S-6 scenarios that contain irrelevant dimensions, we evaluate performances of algorithms 
only on the basis of relevant gene profiles, 
and the adjust rand index is computed by removing irrelevant genes for the (a), (b), (c), (d) and (e) boxes.
Our CSD algorithm allows to select much relevant genes or the genes in the core stuctures. 
We then evaluated in (f) boxes  the agreement between the irrelevant profiles and 
the profiles that are out  of the core structures.
Thus there is one class label reserved to the outer genes of the core structures for the CSD clustering results, 
one class label reserved to the irrelevant genes for the true cluster labels and the rand index is computed on all genes.
The WGCNA algorithm also enables us to exclude irrelevant variables, and as describe above,
we evaluated in (g) boxes the ability of the algorithm to detect irrelevant variables}
\end{minipage}
 \end{figure}

\noindent 
We evaluated the quality of the clusters formed by the three clustering algorithms with the adjust Rand index \cite{milligan1986study},
comparing the clustering results with the true class labels.
For all the simulated datasets with the six scenarios, 
the results are shown in Figure~\ref{simu_RI} as boxplots of the distribution of the adjust rand index.  
From S-1 simulated data, it is seen that all methods were able to recover the underlying cluster structure well, 
when classes are well-defined (high intra-cluster density and low inter-cluster density), 
even if there is many irrelevant genes (S-5).
However, the selection  of the optimal number of clusters for the spectral clustering may still failed, 
especially in the presence of irrelevant genes.
Because of using a single power level to compute weighted adjacency matrix,
the WGCNA algorithm is clearly less adapted to detect clusters with various intra or inter densities 
(S-2, S-3), that is, however, a realistic case. 
Moreover, when cluster may have high inter densities (S-3), 
the choice of optimal number of clusters for spectral clustering is more difficult. 
The three algorithms are able to detect clusters in presence of a reasonable level of noise (S-4), 
but for the WGCNA algorithm it was necessary however, to increase cut height level  to have good results.
Globally, the clustering accuracy is better with our CSD algorithm (S-6), 
which is able to detect classes of different densities and  less sensitive to choice of the parameter in the presence of noise.
The CSD algorithm have also the advantage of allowing a selection of more relevant genes composing the core structures. 
If clusters are very well defined (S-1), the core structures coincide with the clusters (all the genes are in a core structure), 
and even if not all genes are in core structures (S-2, S-3, S-4, S-5),  a minority of them are excluded (values of rand index in (b) boxes remain high). 
The identification of core structures offer the possibility to exclude irrelevant genes, 
while many irrelevant genes are not excluded with the WGCNA algorithm (see (f) and (g) boxes for S-5 and S-6 scenarios).


\subsection{Performance evaluation on real gene expression data}

Using simulated data we have shown that the CSD algorithm has very good performance 
and tends to outperform the spectral and WGCNA algorithms to cluster high network.
Now we are going to look at the performance of CSD algorithm on real gene expression data.
For this purpose we use the well-known synchronized yeast cell cycle data set of Spellman et al. \cite{spellman1998comprehensive},
which include 77 samples under various times during the cell cycle
and a total of 6179 genes, of which 1660 genes are retained for this analysis after preprocessing and filtering.
We deleted genes with missing values more than 20\% or standard deviation less than 0.4 and 
missing values were imputed by the KNN algorithm \cite{troyanskaya2001missing}.
The goal of cluster analysis is to identify 
the correlation patterns in the time series of gene expressions of yeast measured along the cell cycle.
We defined the elements of the weighted adjacency matrix of the graph as 
the absolute values of Spearman rank correlations between gene expression levels.
We conducted sensitivity analysis, used to test and compare algorithm robustness to parameter settings, 
and evaluated clustering results based on internal and external criteria. 

\begin{figure}[!h]
\centering
\begin{minipage}[b]{0.9\linewidth}
\centering
  \includegraphics[scale=0.6]{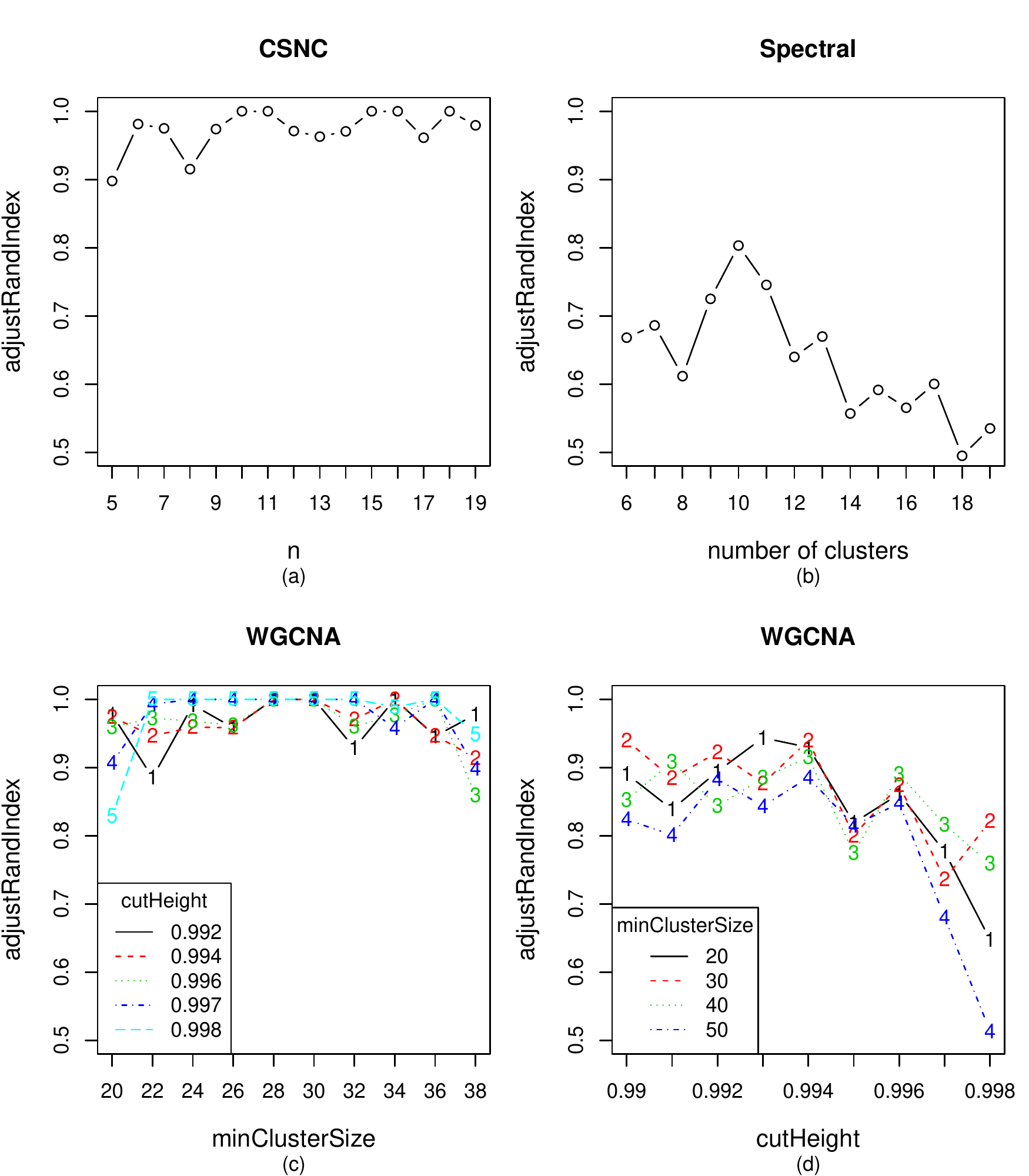}
  \caption{\label{RI_para} \small Parameter sensitivity analysis on yeast cell-cycle gene expression dataset. 
  The adjust rand index was computed to compare the clustering result produced for a given input parameter 
  (on the x-axis) with the consecutive result obtained for an increase in the input parameter.   }
\end{minipage}
\end{figure}

\noindent
To illustrate the robustness of the CSD algorithm to changes of the parameter settings, 
and to compare with the robustness of WGCNA and spectral clustering, 
we computed the adjust rand index between clustering results obtained for consecutive values of the input parameters (Figure~\ref{RI_para}). 
We show that our algorithm is less sensitive to changes in the input parameter and relatively easy to setup, 
while parameter setting is less obvious with the two other algorithms.
Clustering results with spectral clustering are significantly influenced by the clusters number,
and the dynamic tree cut method need to specify two input parameters would affect both 
detection of clusters, should have a wide range of densities, and the number of non-clustered elements. 

\subsubsection{Internal and external evaluations}

We are interested here in evaluating the performance of the clustering algorithms using internal and external criteria. 
Intuitively, based on the data itself (internal clustering), a good clustering will have a high internal density 
(clusters having more edges linking their elements among themselves) and sparser external connections. 
A lot of evaluation metrics have been proposed to quantify the internal quality of a given clustering result.
Most of them are however imperfect because they are typically biased toward one algorithm or the other, 
but may be used as indicative of cluster quality.
We compared algorithms results based on four popular quality metrics, 
the Dunn's index \cite{dunn1973fuzzy}, the modularity \cite{newman2004finding}, the Silhouette width \cite{rousseeuw1987silhouettes}, and the Figure Of Merit (FOM) \cite{yeung2001validating}. 
We will have a more rigorous evaluation of the performance of algorithms by using diversified validation measures. 
We used the adjust Rand index as external validation indices to measure the agreement  between clustering results and  
the biological classes established by Spellman et al. related to five different 
``phase groups'' during the cell cycle

\begin{center}
 \begin{tabular}{|c|c|c|c|c|c|}
  \hline
   \textit{phase group term} & \textbf{G1} & \textbf{G2} & \textbf{M} & \textbf{M/G1} & \textbf{S} \\ \hline
   \textit{group size} & 211 & 63 & 133 & 82 & 42 \\ \hline
 \end{tabular}
 \end{center}

\begin{figure}[!h]
\centering
\begin{minipage}[b]{0.8\linewidth}
\centering
  \includegraphics[scale=0.6]{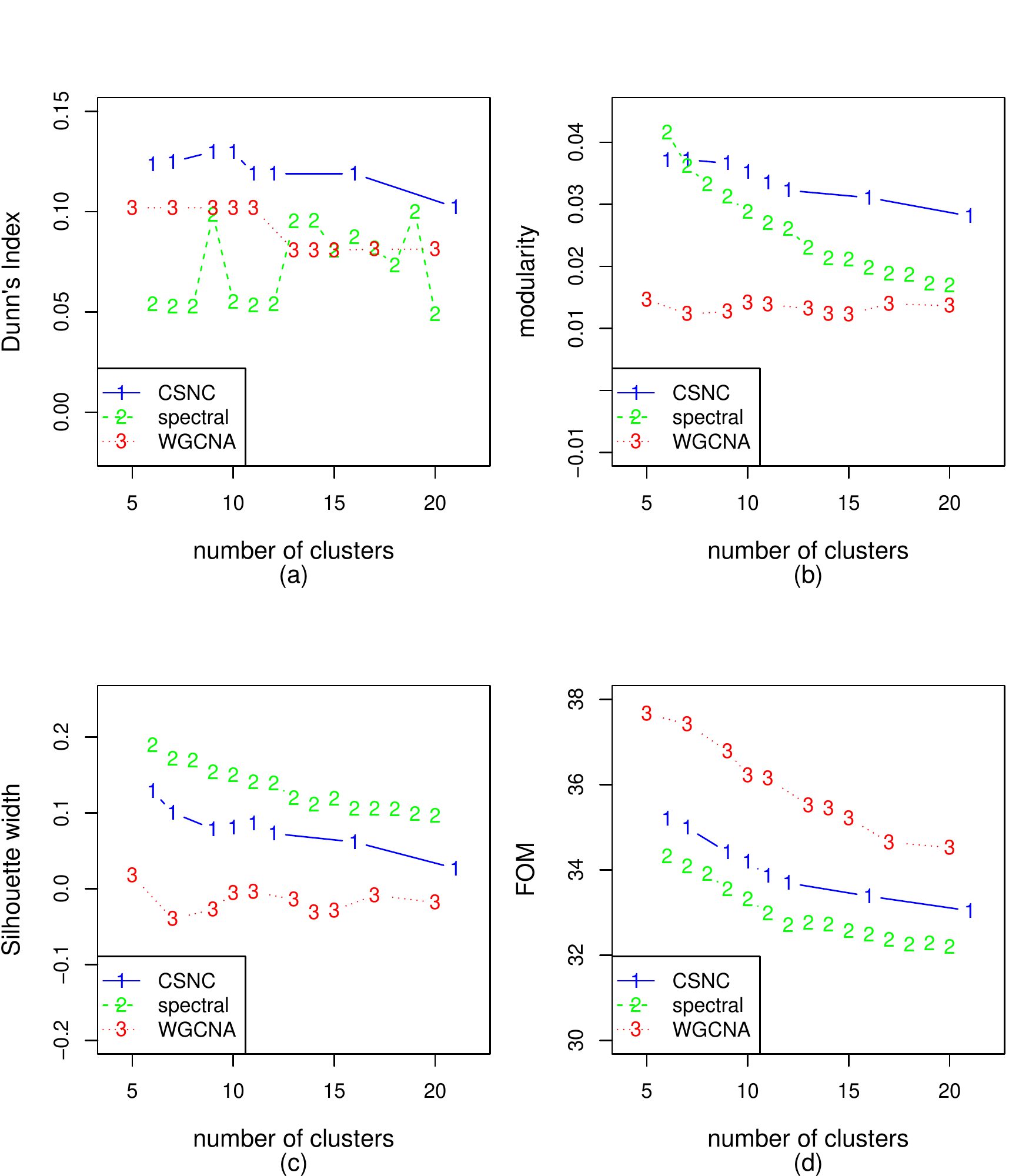}
  \caption{\label{internal_evaluation} \small Internal evaluation of clustering algorithms  on the yeast cell-cycle gene expression data (1622 genes). 
  (a) Dunn's index. (b) Modularity. (c) Silhouette width. (d) Figure Of Merit. }
  \end{minipage}
\end{figure}

\begin{figure}[!h]
\centering
\begin{minipage}[b]{0.7\linewidth}
\centering
  \includegraphics[scale=0.35]{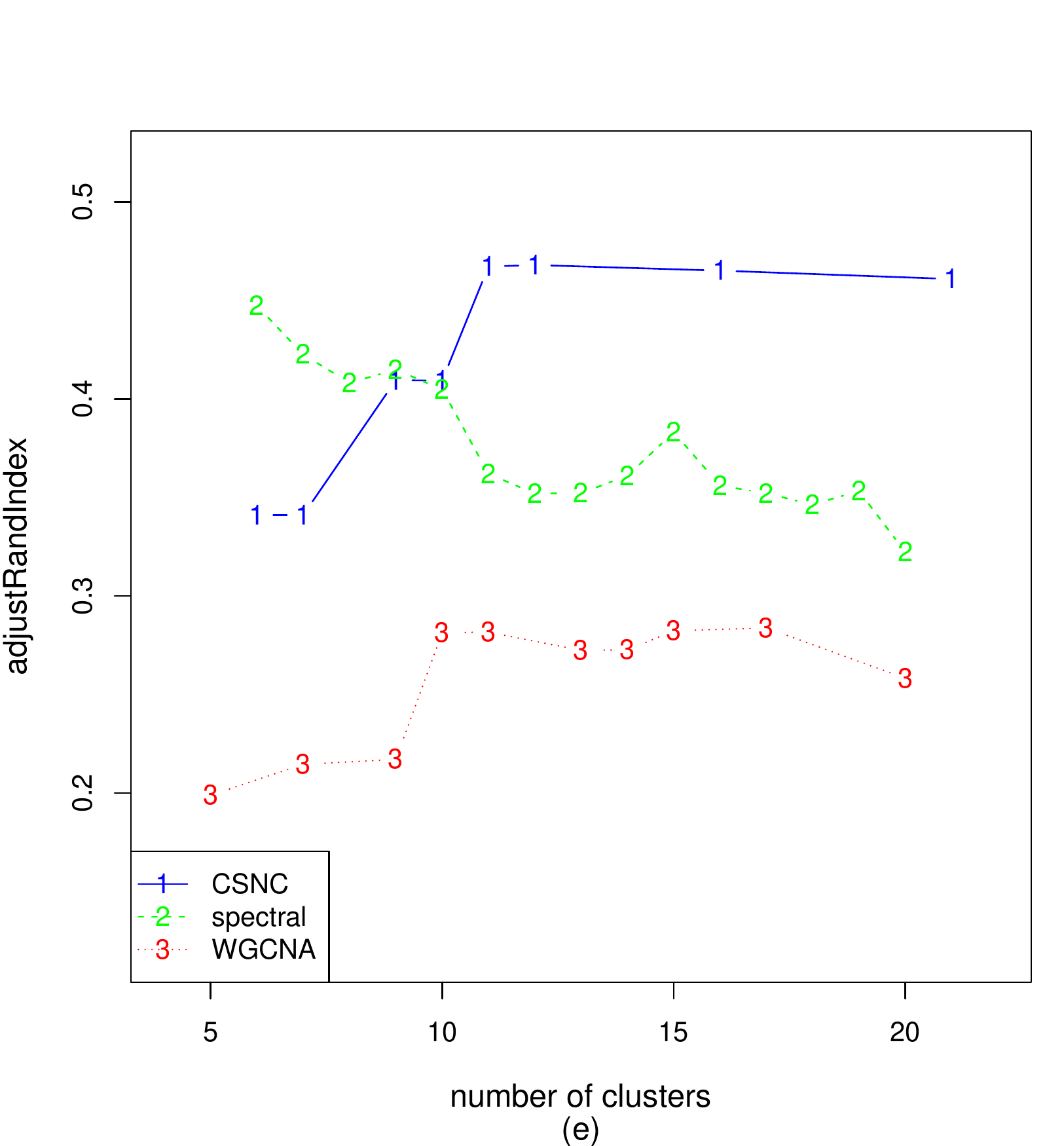}
  \caption{\label{external_evaluation} \small  External evaluation of clustering algorithms  on the yeast cell-cycle gene expression data (531 genes).  
  The adjust rand index to compare clustering results with the gold standard classes established by Spellman.   }
  \end{minipage}
  \end{figure}

\noindent
We compared the results obtained by the three algorithms while varying the number of cluster.
There is a monotonic relationship between the input parameter of the CSD algorithm and the number of cluster formed,
which permits an easy comparison with spectral clustering results.
We used the dynamic tree cut algorithm to cut the hierarchical tree obtained with the WGCNA algorithm.
The height parameter of the dynamic tree cut method has been chosen sufficiently high (0.996 or 0.997) 
to have a maximum number of clustered genes and we controlled the number of clusters by varying the minimum cluster size parameter. 
The subset of genes (1622 genes) clustered with WGCNA algorithm has been used for internal evaluation of algorithm results. 
Out of these genes, 531 are found in gold standard classes established by Spellman and have been used for external clustering validation.
\\ \newline
Our CSD algorithm achieves the best compromise between all internal validation indices (Figure~\ref{internal_evaluation})
and gives the most biologically meaningful gene groups consistent with previous knowledge (Figure~\ref{external_evaluation}).

\newpage
\subsubsection{Functional validation based on gene ontology}
One of the common assertion in expression analysis is that genes sharing a similar pattern of expression 
are more likely to be involved in the same regulatory processes \cite{d2000genetic}. 
Prior knowledge about the biological functions of some genes are available in the 
Gene Ontology (GO) database  which is organized in three levels of detail 
for the decription of the biological process, molecular function, and cellular component of gene products.
We have limited our analysis to the GO terms which describe biological processes. 
Each gene may be associated to one or more biological process GO term. 
We performed enrichment analysis to identify biological processes over-represented in each community of genes  \cite{gene2004gene}. 
The enrichment test is based on the hypergeometric test or, equivalently, the one-tailed Fisher Exact Test, 
indicates if a given GO term appears in a sublist of genes  (genes within a given community) higher than expected by chance \cite{subramanian2005gene}. 
\\ \newline
Table~\ref{enrichment} presents the three more enriched biological process that have a Bonferroni corrected p-value $\leq 0.05$ within each cluster 
produced by the CSD algorithm for a minimum core size parameter equal to 10 (12 clusters). 
We confirmed that the results produced by the proposed method are reliable according to the Spellman et al. results.
The cluster 5 of genes involved in DNA repair contains 70\% of genes in group phase G1 and all the genes group in the CLN2 cluster by Spellman and al.
The small cluster 10 contains all genes placed in the Y cluster by Spellman and al. and 15\% of genes in group phase G1.
Following G1 phase during S phase, nucleosomes are assembled, using histones.
The cluster 6 is composed by genes involved in this phase (contains 81\% of genes in group S and all genes identified by Spellman in the histone cluster). 
The cluster 11 groups genes involved in the tricarboxylic acid cycle, an essential cycle in aerobic growth, and contains 33\% of genes group phase G2.
The cluster 4 contains genes involved in mitosis (71\% of genes in group M). 
We found in this cluster the genes in the Spellman's cluster CBL2 and in the cluster MCM which is induced by CLB2.  
The cluster 12 contains genes involved in cytokinesis and which reach peak in late M or at the M/G1 boundary.
All genes in the Spellman's cluster SIC1 and 30\% of genes in M/G1 are placed in this cluster.

\begin{table}[!h]
\centering
\begin{minipage}[b]{1\linewidth}
\centering
  \includegraphics[trim = 1.9cm 10cm 0 2.5cm, clip=true, scale=0.7]{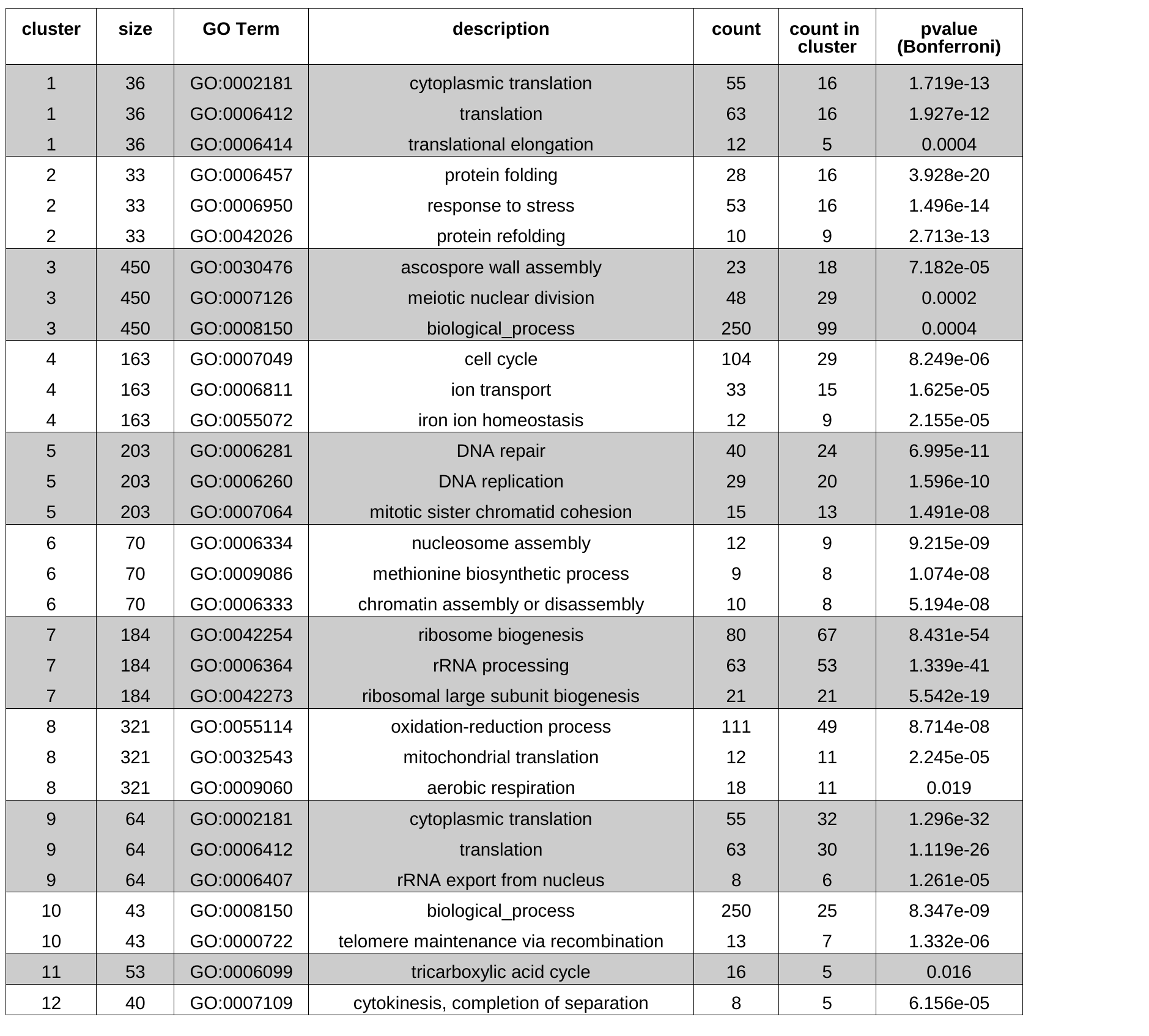}
  \caption{\label{enrichment} \small Enrichment analysis of the 12 clusters produced by the CSD algorithm.   }
\end{minipage}
 \end{table}

\subsection{The core structures include the most biologically informative genes}

We have shown previously that our CSD algorithm is able to identify biologically meaningful gene clusters from a gene co-expression network.
Now, we will show the biological relevance of the core structures. 
It is well accepted that genes with critical functional roles are centrally positioned on the gene co-expression network 
and have high connectivity (hubs).
The core structures by definition contain highly connected or co-expressed genes, 
and so have the best chance of containing potential biologically meaningful genes.
\\ \newline
We studied the core structures constructed by the CSCN algorithm on the yeast cell-cycle gene expression data set of Spellman et al.,
and compared the biological interpretation of the core structures with interpretations based on clustering analysis of all the genes 
presented in the previous section.
Within the core structure we found 472 genes out of a total of 1660 genes.   
We will make comparisons based on two levels of information, enriched GO terms 
and detection of cell cycle-regulated genes identified by Spellman.
We show in Table~\ref{enrichment_core} the three more enriched biological process that have a Bonferroni corrected p-value $\leq 0.05$ 
within each core stucture produced by the CSD algorithm for a minimum core size parameter equal to 10 (12 clusters).
Annotated biological process for core structures (Table~\ref{enrichment_core}) and for entire clusters (Table~\ref{enrichment}) are very similar, 
demonstrating that the core structures hold the most meaningful information.
In addition to that, the core structures highlight most of the genes identified by Spellman to be cell cycle-regulated. 
Among the 472 genes located in the core structures, we found 85\% of the 156 genes 
identified in the previous analysis on the 1660 genes and clustered in the CLN2, Y, histone, CBL2, MCM and SIC1 by Spellman. 
\begin{table}[!h]
\centering
\begin{minipage}[b]{1\linewidth}
\centering
  \includegraphics[trim = 1.8cm 14cm 0 2.5cm, clip=true, scale=0.7]{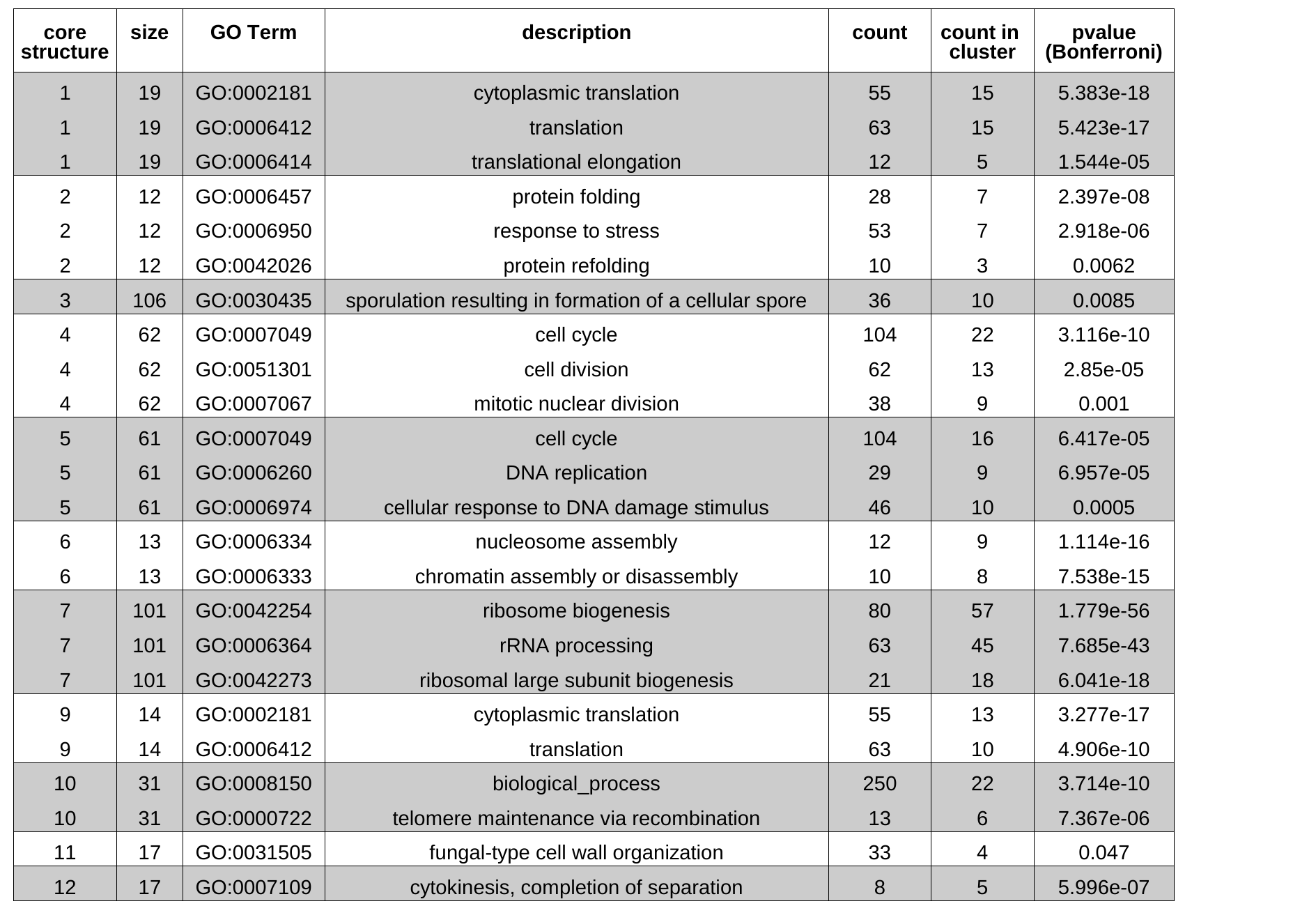}
  \caption{\label{enrichment_core} \small Gene Ontology enrichment analysis for the 12 core structures produced by the CSD algorithm 
  applied on the yeast cell-cycle gene expression data set of Spellman et al. }
  \end{minipage}
\end{table}

\begin{figure}[!h]
\centering
\begin{minipage}[b]{0.9\linewidth}
\centering
  \includegraphics[scale=0.8]{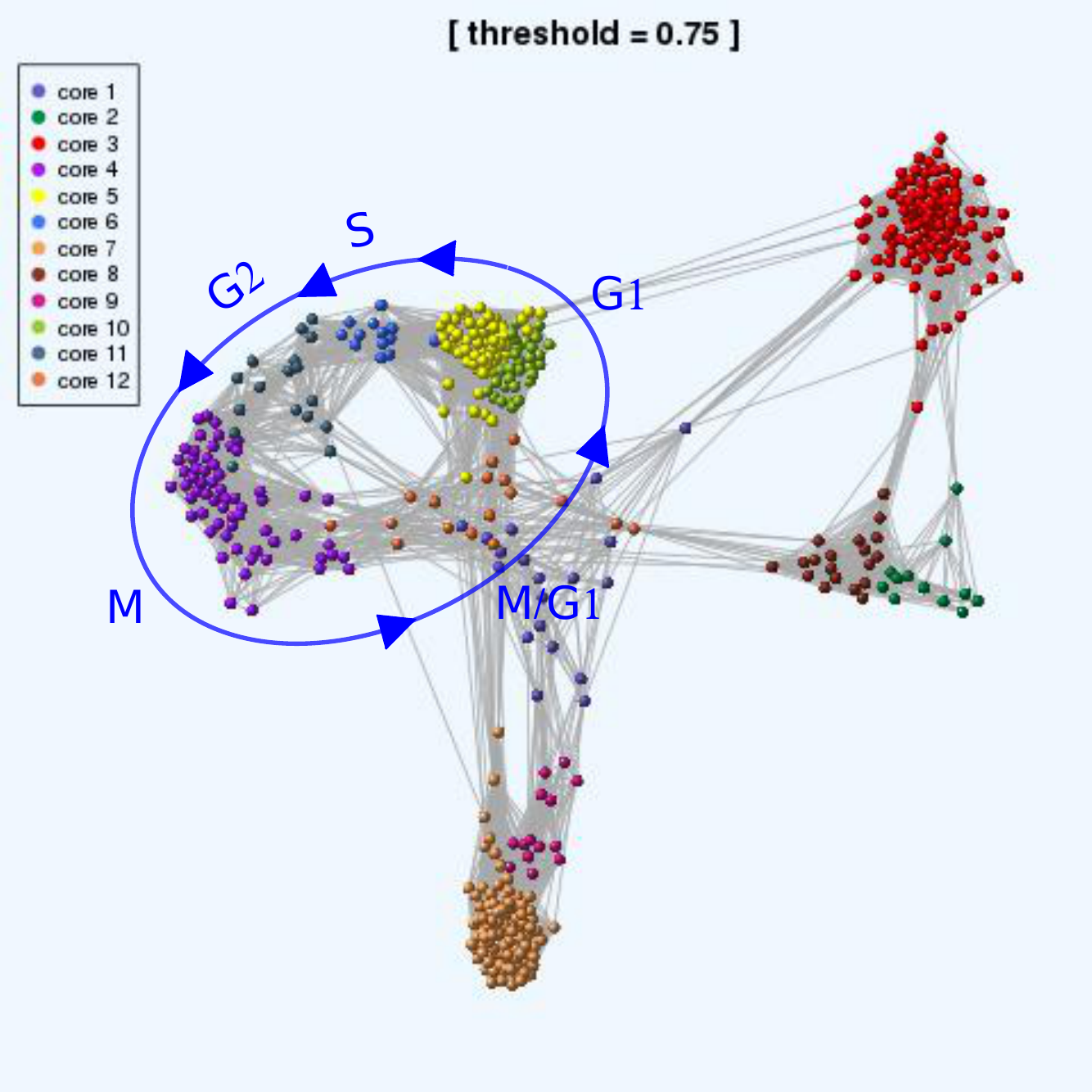}
  \caption{\label{network_core} \small Gene co-expression network of genes included in the core structures produced by the CSD algorithm 
  applied on the yeast cell-cycle gene expression data set of Spellman et al. 
  Each core structure is shown with a different color and an edge connect two nodes if its weight (similarity between the two genes) 
  is greater than 0.75.}
  \end{minipage}
\end{figure}

\noindent
In addition, we have represented on Figure~\ref{network_core} 
the induced subgraph of the graph containing nodes included in the core structures, 
by keeping only edges with a weight larger than 0.75, to have a sparse representation of this subgraph.
This simple representation gives us an idea of the interactions between the different core structures, 
in line with biological knowledge.
Indeed, it is very interesting to observe a ring describing the phases of the cell cycle.
Core structures have been previously associated with one phase of the cell cycle and 
are placed on the ring in the order of the cycle that suggests that genes in a given core structure interact strongly with genes in other core structures 
associated to the same phase or consecutive phases of the cell cyle, which makes sense. 

 
\section{Conclusions}
Investigation of complex interaction patterns among genes is an effective way   
for understanding complex molecular process, predicting gene functions,
finding genes with critical functional roles.
The discovery of clusters of co-expressed genes and centrally positioned genes  is of fundamental and practical interest.
In this paper, we presented a promising algorithm for both the analysis of the community structure in large interaction networks 
and the selection of subgroups of genes centrally positioned in the network. 
Several method for network clustering exist.
However, they have their limitations and require to make assumptions can be restrictive in practice.
For example, the WGCNA algorithm (hierarchical clustering) needs a sparse estimation of the graph 
by thresholding the input similarity matrix (e.g. absolute value of Pearson correlations), 
needs to choose an appropriate similarity measure between nodes of the graph 
and needs to select other parameters for tree cutting. 
These choices can have a large impact on the results.
Another common method for graph clustering is the spectral clustering (partitional clustering).
The main limitation of this method is the difficult problem of determining the number of clusters.
Another drawback of this algorithm is that it tends to form clusters of homogeneous sizes. 
\\ \newline
Our Core Structures based Network Clustering (CSD) algorithm combines advantages of existing methods 
and reduce the bias introduced by limiting the number of input parameter.
Initially, the variables (e.g. gene expressions) are organized in a similarity matrix 
(e.g. absolute value of Pearson correlations). 
The first step of the algorithm is to identify group of variables preferentially connected, called core structures, 
on the collection of all neighborhood graphs obtained by similarity thresholding. 
The algorithm doesn't need to determine an unique threshold parameter for graph estimation by studying all possible graphs.
A unique input parameter is required which is the minimum size of core structures.
This parameter controls the level of granularity of the results (embedded results) 
and this is easy to choose it in order to have not too many small clusters and not too large clusters. 
The second step of the algorithm consists in clustering all the variables based on the core structures knowledge.
Each cluster is composed of a unique core structure and nodes associated with this core structure in the 
iterative process consists of selecting the edge of maximum weight between one clustered element and one unclustered element
and attributing the cluster label of the clustered element to the unclustered element.  
\\ \newline
We evaluated performances of our CSD algorithm and compared with WGCNA and spectral clustering algorithms using simulated data and 
public gene expression data of yeast cell cycle have been studied by Spellman \cite{spellman1998comprehensive}.
We showed that our algorithm outperforms other algorithms based of internal and external criteria, 
is less sensitive to the input parameter, is easy to use 
and is able to select groups of centrally positioned genes in the network or core structures 
contain informative genes.



\section{Methods}

\subsection{Network clustering}

\subsubsection{CSD algorithm}
Let $S=(s_{ij})_{1\leq i,j \leq p}$ be a similarity matrix 
which defines a fully-connected graph $G$ whose weighted adjacency matrix is $W=S$. 
The input parameter or the minimum size of a core structure is denoted by $n$.
Our CSD algorithm is described in algorithm~\ref{CSD_algo}. 
There are three main steps for CSD algorithm,
identification of the maximum spanning tree of G (e.g. Prim's algorithm), 
identification of the core structures (algorithm~\ref{core_algo}) 
and clustering based on the core structures knowledge (algorithm~\ref{clust_algo}).

\begin{algorithm}
\caption{CSD algorithm}\label{CSD_algo}
\begin{algorithmic}[1]
\Procedure{CSD}{W,n} \Comment{\textit{W is the weighted adjacency matrix ; n is the minimum size of a core structure}}
\State $T\leftarrow$ MaximumSpanningTree(W) \Comment{\textit{algorithm (Prim's algorithm) returns the weighted adjacency matrix of the maximum spanning tree}}
\State $Q \leftarrow$ CORE(T,N)
\State $\{C_1,...,C_K\} \leftarrow$ CLUST(Q,T)
\State \textbf{return} $\{C_1,...,C_K\}$\Comment{\textit{The family of clusters}}
\EndProcedure
\end{algorithmic}
\end{algorithm}

\begin{algorithm}
\caption{Core structure algorithm}\label{core_algo}
\begin{algorithmic}[1]
\Procedure{CORE}{W,n} \Comment{\textit{$W\in \mathbb{R}^{p\times p}$ }}
\State $Q \leftarrow \{1,...,p\}$ \Comment{\textit{Q is the family of core structures}}
\While{$|\{ (i,j)\ ;\ i,j\in \{ 1,...,p\},\ w_{ij}>0 \}|> 2(N-1)$}
\State Selecting the lowest cost edge $(i,j)$ of the graph defined by $W$ 
\State $w_{ij}\leftarrow 0$
\State $V_i \leftarrow$ DFS($W,i$) \Comment{\textit{Depth First Search algorithm (DFS) returns nodes in connected component of the graph defined by W containing i}}
\State $V_j \leftarrow$ DFS($W,j$)  
  \If {$|V_i|\geq n \mbox{ and } \ |V_j|\geq n $}
  \State Selecting the core structure $S\subset Q$ containing $V_i\cup V_j$
  \State $Q \leftarrow Q \setminus S$ 
  \State $Q \leftarrow \{Q , V_i , V_j \}$ 
  \Else
  \If {$|V_i|\leq n$}
  \ForAll{$u,v\in V_i$}
    $w_{uv}\leftarrow 0$
  \EndFor
  \EndIf
  \If {$|V_j|\leq n$}
    \ForAll{$u,v\in V_j$}
    $w_{uv}\leftarrow 0$
  \EndFor
  \EndIf
  \EndIf
  
\EndWhile
\State \textbf{return} $Q$\Comment{\textit{The family of core structures}}
\EndProcedure
\end{algorithmic}
\end{algorithm}

\begin{algorithm}
\caption{Clustering algorithm}\label{clust_algo}
\begin{algorithmic}[1]
\Procedure{CLUST}{Q,W} 
\Comment{$Q=\{ Q_k ; k=1,...,K\}$ }
\ForAll {$k = 1,...,K$ }
 \State $C_k \leftarrow Q_k$
\EndFor
\State $U \leftarrow \{i\ ;\ i\in\{1,...,p\},\  i\not\in Q_k ,\forall k=1,...K  \}$
\While{$U\neq \emptyset$}
\State Selecting the greatest cost edge $(i,j)$ of $G$ such that one node is in $U$ and the other node is in $C_k$ ($k=1,...,K$).
Suppose that $i\in U$ and $j\in C_k$.
\State $C_k \leftarrow C_k\cup \{ i\}$
\State $U \leftarrow U\setminus \{i\}$
\EndWhile
\State \textbf{return} $\{C_1,...,C_K\}$\Comment{\textit{The family of clusters}}
\EndProcedure
\end{algorithmic}
\end{algorithm}

\subsubsection{Spectral clustering}

Spectral clustering \cite{ng2002spectral,von2007tutorial} consists in changing the representation of the initial data set 
to enhance the clustering structure and make it easy to detect via standard algorithms, 
like k-means clustering, in the new representation.
The main tool for spectral clustering is the graph Laplacian matrix. 
Let $W=(w_{ij})_{1\leq i,j \leq p}$ be the weighted adjacency matrix of the graph.
The degree matrix $D$ is a diagonal matrix, whose diagonal is the degree for each node, 
the sum of the weights of edges are incident with the node, $d_{ii}=\sum_{j\neq i} w_{ij}$.  
We applied spectral clustering method by using the normalized Laplacian $L_{rw}$ \cite{chung1997spectral} defined as $L_{rw}=D^{-1}L$ 
where $L=D-W$ is the unnormalized Laplacian. 
The change of the representation is induced by the eigenvectors of $L_{rw}$.
In order to partition a graph into K classes,
we compute with the first K eigenvectors (corresponding to the K smallest eigenvalues) of the normalized Laplacian $L_{rw}$ 
and cluster the points in $\mathbb{R}^K$ whose coordinates are elements of eigenvectors, 
with the k-means algorithm into K clusters (R  eigen() and kmeans() functions).

\subsubsection{WGCNA}

All the needed functions for network clustering in WGCNA procedure \cite{zhang2005general} are available in the WGCNA R software package \cite{langfelder2008wgcna}. 
The  weighted adjacency matrix $W=(w_{ij})_{1\leq i,j \leq p}$ of the graph is defined 
by transforming the initial similarity matrix $S=(s_{ij})_{1\leq i,j \leq p}$ using a procedure 
consists in raising the similarity to a power, $w_{ij}=s_{ij}^\beta$. 
To choose the threshold parameters $\beta$ we applied the proposed approximate scale-free topology criterion.
As suggested by Langfelder and Horvath, we performed hierarchical clustering using the topological overlap measure (TOM) as input similarity, 
and used the Dynamic Tree cut algorithm to extract clusters from the dendrogram.

\subsection{Cluster evaluation}
\subsubsection{Internal evaluation}
The Dunn's index \cite{dunn1973fuzzy} the goodness of a clustering, by measuring the maximal diameter of clusters 
and relating it to the minimal distance between clusters.
The Dunn's index $D_K$ for a given clustering partition, $C_1,...,C_K$ is defined as follows
\begin{equation} \label{dunn_eq}
 D_K = \frac{\displaystyle\min_{1\leq i < j \leq K}\ \delta (C_i,C_j)}{\displaystyle\max_{1\leq k\leq K}\ \Delta_k},
\end{equation}
where $K$ is the number of clusters, $\delta (C_i,C_j)$ is the intercluster dissimilarity between $C_i$ and $C_j$ defined as 
$\delta (C_i,C_j)=\displaystyle\min_{x\in C_i, y\in C_j} d(x,y)$, and $\Delta_k$ is the diameter of the cluster $C_k$ defined as
$\Delta_k=\displaystyle\max_{x,y\in C_k} d(x,y) $.
Thus, high Dunn's index values indicate the presence of clusters such that the distances between the clusters are high and 
the diameter of the clusters are small.
The Dunn's index does not exhibit any trend with respect to number of clusters and can be used as 
an indication to help the choice of the number of clusters for a clustering algorithm needs to tune this parameter.
\\ \newline
The Silhouette width \cite{rousseeuw1987silhouettes} 
is another index that combine the compactness and the separation measures. 
A silhouette value is measured for each node of the graph to quantify
the degree of confidence in the clustering assignment of the nodes. 
For one node $i$, the silhouette is defined as follows
\begin{equation} \label{silhouette_eq}
S(i)=\frac{b_i - a_i}{\max(b_i,a_i)},
\end{equation}
where $a_i$ is the average dissimilarity between $i$ and all the other nodes in the same cluster,
and $b_i$ is the average dissimilarity between $i$ and the observations in the nearest cluster 
$$
b_i=\displaystyle\min_{ \{ C_k ; i\not\in C_k \} }\displaystyle\sum_{j\in C_k} \frac{d(i,j)}{|C_k|},
$$
where $d(i,j)$ is the dissimilarity between nodes i and j.
The Silhouette width $S$ is then the average value of each node's silhouette value, should be maximized and 
lies in the interval $[-1,1]$.
\\ \newline
A figure of merit (FOM) \cite{yeung2001validating} is an estimate of the predictive power of a clustering algorithm.
A typical gene expression data set contains measurements of expression levels of p genes in n samples.
The gene expression data matrix is denoted by $X\in \mathbb{R}^{n\times p}$. 
Suppose there are K clusters, $C_1,C_2,...,C_K$. 
Let $x_{ij}$ be the expression level of the gene $j$ in the sample $i$, 
and $\mu_i^{C_k}=\frac{1}{|C_k|}\sum_{j\in C_k}x_{ij}$ be the average expression level in the sample $i$ of genes in cluster $C_k$.
This mean expression $\mu_i^{C_k}$ is used as if it was a prediction of expression level in the sample $i$ 
for each gene in the cluster $C_k$ and is compared with the known values of expression in this sample for all these genes.
The figure of merit $FOM$ is defined as follows
\begin{equation} \label{FOM_eq}
FOM=\displaystyle\sum_{i=1}^n\sqrt{ \frac{1}{p} \displaystyle\sum_{k=1}^K \displaystyle\sum_{j \in C_k} (x_{ij}-\mu_i^{C_k})^2 }. 
\end{equation}
The lower the FOM value is, the higher the predictive power of the algorithm is.
\\ \newline 
We used the dissimilarity measure $d(i,j)=1-x_{ij}$ for both the dunn's index and the silhouette width, 
where $w_{ij}\in[0,1]$ is element of the weighted adjacency matrix (absolute values of Pearson or Spearman correlations).
The Dunn's index, the silhouette width and the FOM are implemented in the R package clValid \cite{brock2011clvalid}.
\\ \newline
The modularity of Newman and Girvan \cite{newman2004finding} is the most popular quality function for graph clustering 
that measures the non-randomness of a graph partition. 
This index has been used in many algorithms as a quality function or as an objective function \cite{newman2004fast,ruan2010general,ruan2007efficient,white2005spectral}. 
The idea behind the modularity index is that a graph have a communitity structure if 
it is different from a random graph or a null model is not expected to have communitity structure.
The null model proposed by Newman and Girvan consists of a randomized version of the original graph, 
where edges are placed at random, under the constraint that the expected degree of each node matches 
the degree of the node in the original graph.
The modularity index $Q$ is defined as follows
\begin{equation} \label{modularity_eq}
Q=\displaystyle\sum_{i=1}^K \left[ \frac{l_i}{m} - \left( \frac{d_i}{2m} \right)^2 \right],
\end{equation}
where $K$ is the number of clusters, 
and for the weighted version of the modularity, 
$m$ is the sum of weights of all edges of the graph, 
$l_i$ is the sum of weights of edges connecting vertices of module $C_i$ 
and $d_i$ is the sum of the degrees of the vertices of $C_i$.
A high modularity means that there are more edges within the clusters that you expect by chance.
Modularity is always smaller than one and can be negative.

\subsubsection{External evaluation}
To compare the clusters identified by an algorithm to the gold standard clusters, 
we computed the adjusted Rand index \cite{milligan1986study}.
Suppose that $A_1,A_2,...,A_{K_A}$ and $B_1,B_2,...,B_{K_B}$ represent two different partitions of the nodes of the graph. 
Let $n=|V|$ be the number of nodes in the graph and $n_{ij}$ be the number of nodes that are both in cluster $A_i$ and in cluster $B_j$. 
Let $n_{i.}=|A_i|$ and $n_{.j}=|B_j|$ be the number of nodes in cluster $A_i$ and in cluster $B_j$ respectively.
The adjusted Rand index is defined as follows
\begin{equation} \label{RI_eq}
 \frac{\sum_{i,j} {n_{ij} \choose 2} - \left[ \sum_i {n_{i.} \choose 2}\sum_j {n_{.j} \choose 2} \right] / {n \choose 2} }
 { \frac{1}{2}\left[ \sum_i {n_{i.} \choose 2} + \sum_j {n_{.j} \choose 2}\right] 
 - \left[ \sum_i {n_{i.} \choose 2} \sum_j {n_{.j} \choose 2} \right] / {n \choose 2}  }
\end{equation}
The adjusted Rand index lies in the interval $[0,1]$ and is equal to $1$ when the two partitions agree perfectly.

\subsubsection{Functional evaluation}
To assess the functional significance of gene clusters, 
we analysed the enrichment of GO terms for the genes within each cluster. 
We used the R packages org.Sc.sgd.db \cite{package_orgCssgddb} and GO.db \cite{package_godb} R package in Bioconductor to identify GO terms associated with the genes. 
We measured the statistical significance of the GO term enrichment by an hypergeometric test (R phyper() function).
The p-values are adjusted by Bonferroni corrections for multiple testing problems.

\bibliographystyle{plain}
\bibliography{article}

\end{document}